\documentclass[11pt]{amsart}
\usepackage{amscd,amssymb,verbatim}
\theoremstyle{plain}
\newtheorem{Thm}{Theorem}[section]
\newtheorem{Cor}[Thm]{Corollary}

\newtheorem{Prop}[Thm]{Proposition}

\newtheorem{remark}{Remark}

\begin{document}
\title[On contractively complemented subspaces]
{On contractively complemented subspaces of 
separable \(L_1\)-preduals}
\author{Ioannis Gasparis}
\address{Department of Mathematics \\
         Oklahoma State University \\
         Stillwater, OK 74078-1058 \\
         U.S.A.}
\email{ioagaspa@math.okstate.edu}
\begin{abstract}
It is shown that for an \(L_1\)-predual space \(X\) and a countable
linearly independent subset of ext(\(B_{X^*}\)) whose norm-closed
linear span \(Y\) in \(X^*\) is \(w^*\)-closed, there exists a 
\(w^*\)-continuous contractive projection from \(X^*\) onto \(Y\).
This result combined with those of Pelczynski and Bourgain yields
a simple proof of the Lazar-Lindenstrauss theorem that every separable
\(L_1\)-predual with non-separable dual contains a contractively
complemented subspace isometric to \(C(\Delta)\), the Banach space
of functions continuous on the Cantor discontinuum \(\Delta\).

It is further shown that if \(X^*\) is isometric to \(\ell_1\)
and \((e_n^*)\) is a basis for \(X^*\) isometrically equivalent
to the usual \(\ell_1\)-basis, then there exists a \(w^*\)-convergent
subsequence \((e_{m_n}^*)\) of \((e_n^*)\) such that the closed
linear subspace of \(X^*\) generated by the sequence
\((e_{m_{2n}}^* - e_{m_{2n-1}}^*)\) is the range of a \(w^*\)-continuous
contractive projection in \(X^*\). This yields a new proof of Zippin's
result that \(c_0\) is isometric to a contractively complemented subspace
of \(X\). 
\end{abstract}
\maketitle
\section{Introduction} \label{S:1}
A Banach space \(X\) is said to be an \(L_1\)-predual
provided its dual \(X^*\) is isometric to \(L_1(\mu)\)
for some measure space \((\Omega, \Sigma, \mu)\).
Perhaps the most natural example of an \(L_1\)-predual
is \(C(K)\), the Banach space of real-valued functions 
continuous on the compact Hausdorff space \(K\), under
the supremum norm. \(L_1\)-preduals were the subject of
an extensive study in the late 1960's and early 1970's.
For a detailed survey of results on \(L_1\)-preduals
we refer to \cite{La}. For the connection between 
\(L_1\)-preduals and infinite-dimensional convexity
we refer to the recent survey article \cite{FL}. 

For some time it was thought that every \(L_1\)-predual
is isomorphic to a \(C(K)\)-space for suitable \(K\), but
the example given by Benyamini and Lindenstrauss \cite{Be}
disproved this. The present paper is concerned with the existence
of subspaces of a separable infinite-dimensional \(L_1\)-predual 
\(X\), isometric to \(C(K)\)-spaces. It was proven by Zippin
\cite{Z} that \(X\) contains a contractively complemented subspace
isometric to \(c_0\). When \(X^*\) is non-separable, Lazar and
Lindenstrauss \cite{LL} proved that \(X\) contains a contractively
complemented subspace isometric to \(C(\Delta)\), where \(\Delta\)
denotes the Cantor discontinuum. These results complement each
other in the sense that neither of them implies the other.  

In the present paper we demonstrate a unified approach towards
these results. Our method consists of establishing Theorem
\ref{T:2} which describes a technique for constructing
a strictly increasing sequence \((V_n)\) of finite-dimensional
subspaces of \(X^*\), where each \(V_n\) is isometric to some
\(\ell_1^{k_n}\), for which there exists an almost 
commuting sequence \((P_n)\)
of \(w^*\)-continuous contractive projections in \(X^*\) such that
\(\mathrm{Im} \, P_n = V_n\), \(n \in \mathbb{N}\).
The proof of Theorem \ref{T:2} is an elementary application of
the principle of local reflexivity \cite{LR}, \cite{JR}.
We apply Theorem \ref{T:2} in order to provide an alternative
proof for the following result
\begin{Thm} \label{T:1}
Let \(X\) be an \(L_1\)-predual and let \(K\) be a countable
subset of \(\mathrm{ext}(B_{X^*})\) such that \(K \cap (-K) = \emptyset\).
Suppose that the norm-closed linear span \(Y\) of \(K\)  
is \(w^*\)-closed in \(X^*\). Then there exists a \(w^*\)-continuous
contractive projection from \(X^*\) onto \(Y\).
\end{Thm}
We remark that Theorem \ref{T:1} is a consequence of the 
following result
formulated by Lazar and Lindenstrauss as Corollary 1
in \cite{LL}:

{\em Suppose that\/} \(X\) {\em is a Banach space so that\/} 
\(X^{*}\)
{\em is isometric to\/} \(L_1(\mu)\). {\em Let\/} 
\(F\) {\em be a face of\/}
\(B_{X^{*}}\) {\em and denote by\/} \(H\) {\em the convex hull of\/}
\(F \cup (-F)\). {\em Assume that\/}
\(H\) {\em is\/} \(w^{*}\)-{\em closed
and metrizable. Then, there exists a\/} \(w^{*}\)-{\em continuous,
affine, symmetric retraction of\/}
\(B_{X^{*}}\) {\em onto\/} \(H\).

Evidently, this result yields Theorem \ref{T:1}.
However, the authors of \cite{LL} offer no proof of this result
and moreover, via their preceding discussion, seem to require that
the face \(F\) be \(w^{*}\)-closed.
Note that if \(F\) is \(w^{*}\)-closed, so is \(H\) but the
converse is not true in general.
We also note that the proof of the aforementioned result that appears
in \cite{La} is false. Specifically, the map \(\phi\)
defined in the proof of the Corollary on page 224 of \cite{La}
is not convex.

Theorem \ref{T:1} has applications in the study of \(\ell_1\)-preduals
\cite{A1}, \cite{G}.
As a consequence of Theorem \ref{T:1} we obtain
\begin{Cor} \label{C:1}
Let \(X\) be a separable \(L_1\)-predual and let \(K\)
be a countable, \(w^*\)-compact subset of 
\(\mathrm{ext}(B_{X^*})\) such that \(K \cap (-K) = \emptyset\).
Then there exists a contractively complemented subspace of \(X\)
isometric to \(C(K)\).
\end{Cor} 
This corollary combined with the results of Pelczynski \cite{P2}
and Bourgain \cite{Bo} yields the next
\begin{Cor} \label{C:2}
Let \(X\) be a separable \(L_1\)-predual such that \(X^*\) is non-
separable. Then there exists a contractively complemented subspace
of \(X\) isometric to \(C(\Delta)\).
\end{Cor}
This result was obtained in \cite{LL} with a different method.
Their proof is based on a remarkable affine version of Michael's selection
theorem \cite{M}
and makes use of a non-trivial result established in \cite{Lz}.

Another application of Theorem \ref{T:2} is the following
\begin{Thm} \label{T:3}
Let \(X\) be a Banach space such that \(X^*\) is isometric to
\(\ell_1\), and let \((e_n^*)\) be a basis for \(X^*\)  
isometrically equivalent to the usual \(\ell_1\)-basis.
Then there exists a \(w^*\)-convergent subsequence 
\((e_{m_n}^*)\) of \((e_n^*)\) such that the subspace
generated by the sequence \((e_{m_{2n}}^* - e_{m_{2n-1}}^*)\)
is the range of a \(w^*\)-continuous contractive projection in
\(X^*\).
\end{Thm}
Theorem \ref{T:3} combined with Corollary \ref{C:2}
yields an alternative proof of Zippin's result \cite{Z}
\begin{Cor} \label{C:3}
Let \(X\) be a separable infinite-dimensional \(L_1\)-predual.
Then there exists a contractively complemented subspace
of \(X\) isometric to \(c_0\).
\end{Cor}
\section{Preliminaries} \label{S:2}
We shall make use of standard Banach space facts and terminology
as may be found in \cite{LT}. In this section we review some
of the necessary concepts. All Banach spaces under consideration 
will be over the field of real numbers.
By the term {\em subspace} of a Banach space 
\(X\) we shall mean a closed linear subspace. We let
\(B_X\) stand for the closed unit ball of \(X\), while 
\(X^*\) denotes its topological dual. A subspace \(Y\) of \(X\)
is said to be {\em complemented} if it is the range of a 
bounded linear projection \(P \colon X \to X\). When
\(\|P\|=1\), \(Y\) is a {\em contractively} complemented subspace
of \(X\).

\(\ell_1\) denotes the Banach space of absolutely summable
sequences under the norm given by the sum of the absolute values 
of the coordinates. The usual \(\ell_1\)-basis is the Schauder basis
of \(\ell_1\) consisting of sequences having exactly one coordinate
equal to \(1\) and vanishing at the rest of the coordinates.
\(\ell_1^k\), where \(k \in \mathbb{N}\), is the \(k\)-dimensional
subspace of \(\ell_1\) spanned by the first \(k\) members of the
usual \(\ell_1\)-basis. A sequence \((x_n)\) in some Banach space
is {\em isometrically equivalent} to the usual \(\ell_1\)-basis, if
\(\|\sum_{i=1}^n a_i x_i\|= \sum_{i=1}^n |a_i|\), for all
\(n \in \mathbb{N}\) and scalar sequences \((a_i)_{i=1}^n\). 
A finite sequence \((x_i)_{i=1}^k\) in some Banach space is
isometrically equivalent to the usual \(\ell_1^k\)-basis, if
\(\|\sum_{i=1}^k a_i x_i\|= \sum_{i=1}^k |a_i|\), for all
scalar sequences \((a_i)_{i=1}^k\). 

\(c_0\) stands for the Banach space of null sequences under the
norm given by the supremum of the absolute values of the coordinates.
\(\ell_\infty^n\) denotes the Banach space \(\mathbb{R}^n\)
under the norm given by the maximum of the absolute values of the
coordinates.

Given a measure space \((\Omega, \Sigma, \mu)\)
with \(\mu\) positive, \(L_1(\mu)\)
denotes the Banach space of equivalence classes of absolutely
integrable functions on \(\Omega\) under the norm
\(\|f\|=\int_{\Omega} |f| d \mu\).  
\(L_\infty(\mu)\) denotes the Banach space of equivalence classes
of essentially bounded \(\Sigma\)-measurable functions on
\(\Omega\) under the norm 
\(\mathrm{ess \, sup_{\omega \in \Omega}} |f(\omega)|\).

An \(L_1\)-{\em predual} is a Banach space \(X\) such that
\(X^*\) is isometric to \(L_1(\mu)\) for some 
measure space \((\Omega, \Sigma, \mu)\).
According to a result of Pelczynski \cite{P1}, Proposition 1.3,
there exists another measure space
\((\Omega', \Sigma', \nu)\) with \(L_1(\nu)\) isometric
to \(L_1(\mu)\) and such that \(L_1(\nu)^*\) is canonically
isometric to \(L_\infty(\nu)\).
Thus in the sequel, by an \(L_1\)-predual we shall mean a Banach space \(X\)
with \(X^*\) isometric to some \(L_1(\mu)\) such that
\(L_1(\mu)^*\) is canonically isometric to \(L_\infty(\mu)\).

Given a linear topological space \(V\) and
\(A \subset V\), we let
\(\mathrm{co}(A)\) denote the convex hull of \(A\).
Let now \(K\) be a convex subset of \(V\). 
A point \(x \in K\) is called an {\em extreme point},
if whenever \(y\), \(z\) are in \(K\) and 
\(x=ay + (1-a)z\) for some \(0 < a < 1\), then
\(x=y=z\). We let \(\mathrm{ext}(K)\) denote the set
of extreme points of \(K\). 
It is well known that for an \(L_1\)-predual space \(X\)
with \(X^*\) isometric to \(L_1(\mu)\), \(\mathrm{ext}(B_{X^*})\)
consists precisely of functions of the form
\(\sigma \chi_A / \mu(A)\), where \(\sigma \in \{-1,1\}\),
\(A\) is an atom with \(0 < \mu(A) < \infty\) and 
\(\chi_A\) stands for the indicator function of \(A\).

We next recall the important principle of local reflexivity
\cite{LR}, \cite{JR} (cf. also \cite{W}).
\begin{Thm} \label{lr}
Let \(X\) be a Banach space and let \(E \subset X^{**}\)
and \(F \subset X^*\) be finite-dimensional subspaces.
Given \(\epsilon > 0\), there exists an invertible linear
operator \(T \colon E \to X\) such that
\(\|T\| \|T^{-1}|TE\| \leq 1 + \epsilon\),
\(T|E \cap X = id\), and 
\(f(Te)=e(f)\) for all \(f \in F\) and \(e \in E\).
\end{Thm}
\section{A construction of \(w^*\)-continuous contractive projections}
This section is devoted to the proof of Theorem \ref{T:2}
which provides a method of constructing 
\(w^*\)-continuous contractive projections onto certain
finite-dimensional subspaces of \(X^*=L_1(\mu)\),
isometric to \(\ell_1^k\). Repeated applications of 
Theorem \ref{T:2} will in turn
enable us to construct a sequence \((P_n)\) of almost commuting
\(w^*\)-continuous contractive
projections in \(X^*\) such that \((\mathrm{Im} \, P_n)\)
is strictly increasing and \(\mathrm{Im} \, P_n\) is isometric
to some \(\ell_1^{k_n}\), for all \(n \in \mathbb{N}\).
In order to construct \(w^*\)-continuous
projections onto subspaces of \(X^*\) isometric to \(\ell_1\),
we shall make use of the following
\begin{Prop} \label{P:1}
Let \(X\) be a Banach space and let \(Y\) be a \(w^*\)-closed
subspace of \(X^*\). Assume that there exists a net
\(\{Y_\lambda\}_{\lambda \in \Lambda}\) of \(w^*\)-closed
subspaces of \(Y\) with \(Y_{\lambda_1} \subset Y_{\lambda_2}\)
whenever \(\lambda_1 \leq \lambda_2\) in \(\Lambda\), and such that
\(\cup_{\lambda \in \Lambda} Y_\lambda\) is norm-dense in \(Y\).
Assume further that each \(Y_\lambda\) is the range of a \(w^*\)-
continuous projection \(P_\lambda\) in \(X^*\), so that
\(\sup_\lambda \|P_\lambda\| \leq M < \infty\), and
\(\lim_\lambda \sup_{\lambda \leq \mu} \|P_\lambda P_\mu 
-P_\lambda\|=0\). Then there exists a \(w^*\)-continuous
projection \(P\) from \(X^*\) onto \(Y\) with \(\|P\| \leq M\).
\end{Prop}  
\begin{proof}
\(B_Y\) is \(w^*\)-compact. By Tychonoff's theorem we infer that
\(K= \prod_{x^* \in B_{X^*}} M B_Y\) is compact when
endowed with the cartesian topology.
We can now identify \(\{P_\lambda\}_{\lambda \in \Lambda}\)
with a net in \(K\) to obtain a sub-net
\(\{P_{\lambda'}\}_{\lambda' \in \Lambda'}\) of 
\(\{P_\lambda\}_{\lambda \in \Lambda}\) which converges to
an element \(P\) of \(K\). Since \(Y\) is \(w^*\)-closed
in \(X^*\), \(P\) induces a bounded linear operator from 
\(X^*\) into \(Y\), which we still denote by \(P\).
Clearly \(Px^*= w^*-\lim_{\lambda' \in \Lambda'} P_{\lambda'}x^*\),
for all \(x^* \in X^*\) and thus \(\|P\| \leq M\).
Our assumptions yield that \(Px^*=x^*\), for all
\(x^* \in \cup_{\lambda \in \Lambda} Y_\lambda\) and hence
\(P\) is a projection onto \(Y\).

We next demonstrate that \(P\) is \(w^*\)-continuous.
By a classical result \cite{KS} it suffices to show that for
every net \((x_\nu^*)\) in \(B_{X^*}\) such that
\(w^*-\lim_\nu x_\nu^*=\mathbf{0}\), we have that
\(w^*-\lim_\nu Px_\nu^*=\mathbf{0}\). Note that 
\(\|Px_\nu^*\| \leq M\), for all \(\nu\), and let
\(y^* \in MB_Y\) be any \(w^*\)-cluster point of
\((Px_\nu^*)_\nu\). We will show that \(y^*=\mathbf{0}\).
To this end set \(\delta_\lambda = \sup_{\lambda \leq \mu}
\|P_\lambda P_\mu - P_\lambda\|\). Then
\(\|P_\lambda P - P_\lambda\| \leq \delta_\lambda\),
as \(P_\lambda\) is \(w^*\)-continuous. It follows
that \(\|P_\lambda Px_\nu^* - P_\lambda x_\nu^*\| 
\leq \delta_\lambda\), for all \(\lambda\) and \(\nu\),
and thus as \(P_\lambda\) is \(w^*\)-continuous, we obtain
that \(\|P_\lambda y^*\| \leq \delta_\lambda\), for all 
\(\lambda \in \Lambda\). Hence \(Py^*=\mathbf{0}\).
Because \(y^* \in Y\) and \(P\) is a projection onto \(Y\),
we deduce that \(y^*=\mathbf{0}\), completing the proof of the assertion.
\end{proof}
We next pass to the key result which is related to
Lemma 3.1 and Corollary 3.2 of \cite{JR}.
\begin{Thm} \label{T:2}
Let \(X\) be an \(L_1\)-predual and let \(V\) be
a subspace of \(X^*\) isometric to \(\ell_1^k\).
Let \((\delta_i)_{i=1}^n\) be a finite sequence of
positive scalars and assume that there exist \(w^*\)-
continuous linear operators \(T_i \colon X^* \to V\),
\(\|T_i\| \leq 1\), \(i \leq n\), as well as linear
operators \(R_i \colon V \to V\), \(\|R_i\| \leq 1\),
\(i \leq n\), so that \(\|R_i T_n - T_i\| < \delta_i\)
for all \(i \leq n\). Assume further that there exist
collections of vectors \((f_j)_{j=1}^q \subset 
\mathrm{ext}(B_{X^*})\) and \((v_j)_{j=1}^q \subset B_V\),
with \((f_j)_{j=1}^q \) linearly independent, such that
\(\|R_i v_j - T_i f_j \| < \delta_i\) for all
\(i \leq n\) and \(j \leq q\). Then there exists a 
\(w^*\)-continuous linear operator \(T \colon X^* \to V\),
\(\|T\| \leq 1\), such that \(\|R_i T - T_i\| < \delta_i\)
for all \(i \leq n\), and \(Tf_j = v_j\) for all \(j \leq q\).
\end{Thm}
The proof of this result will follow after establishing the next
\begin{Prop} \label{P:2}
Under the hypothesis of Theorem \ref{T:2}, for every \(\epsilon > 0\)
there exists a \(w^*\)-continuous linear operator
\(S \colon X^* \to V\), \(\|S\| \leq 1\), such that
\(\|R_i S - T_i \| < \delta_i + \epsilon\) and
\(\|Sf_j - v_j\| < \epsilon\) for all \(i \leq n\) and 
\(j \leq q\).
\end{Prop}
\begin{proof}
Let \((e_l)_{l=1}^k\) be a basis for \(V\) isometrically
equivalent to the usual \(\ell_1^k\)-basis.
Since \(T_i\) is \(w^*\)-continuous and
\(\|T_i\| \leq 1\), there exist vectors
\((z_{il})_{l=1}^k\) in \(B_X\) such that
\(T_i x^* = \sum_{l=1}^k x^*(z_{il}) e_l\) for all
\(x^* \in X^*\). There also exist scalars
\((a_{ils})\), \(i \leq n\), \(l \leq k\), \(s \leq k\),
such that \(R_i e_l = \sum_{s=1}^k a_{ils} e_s\)
for all \(i \leq n\) and \(l \leq k\). Finally, there
exist scalars \((v_{jl})\), \(j \leq q\), \(l \leq k\),
such that \(v_j = \sum_{l=1}^k v_{jl} e_l\), \(j \leq q\).
Observe that for \(x^* \in X^*\) and \(i \leq n\) we have
\[ \sum_{s=1}^k \biggl ( \sum_{l=1}^k a_{ils}
   x^*(z_{nl}) \biggr ) e_s = \sum_{l=1}^k x^*(z_{nl})
   \sum_{s=1}^k a_{ils} e_s = R_i T_n x^* .\]
Hence, \(\sum_{s=1}^k \bigl | 
\bigl (\sum_{l=1}^k a_{ils} x^*(z_{nl}) \bigr ) -
x^*(z_{is}) \bigr | = \|R_i T_n x^* - T_i x^* \| <
\delta_i\), for all \(x^* \in B_{X^*}\) and \(i \leq n\).
Thus
\begin{equation} \label{1}
\biggl \| \sum_{s=1}^k \rho_s \biggl [ \biggl (
\sum_{l=1}^k a_{ils} z_{nl} \biggr )  - z_{is} \biggr ]
\biggr \| < \delta_i, 
\end{equation}
for all choices of signs
\((\rho_s)_{s=1}^k\) and all \(i \leq n\).

Similarly, \(R_i v_j = \sum_{s=1}^k ( \sum_{l=1}^k a_{ils}
v_{jl}) e_s\) for all \(i \leq n\), \(j \leq q\), and thus
\begin{equation} \label{2}
\sum_{s=1}^k \biggl | \biggl ( \sum_{l=1}^k a_{ils} v_{jl} \biggr )-
f_j (z_{is}) \biggr | = \|R_i v_j - T_i f_j \| < \delta_i,
j \leq q, i \leq n.
\end{equation}
We first show that there exist vectors
\((x_l^{**})_{l=1}^k\) in \(B_{X^{**}}=
B_{L_\infty (\mu)}\) such that
\(x_l^{**}(f_j)= v_{jl}\) for all \(l \leq k\), \(j \leq q\),
and \(\bigl \| \sum_{s=1}^k \rho_s 
\bigl [ \bigl (\sum_{l=1}^k a_{ils} x_l^{**} \bigr )
 - z_{is} \bigr ] \bigr \| \leq \delta_i\), for all \( i \leq n\)
and all choices of signs \((\rho_s)_{s=1}^k\).
Indeed, let \(\sigma_1, \dots, \sigma_q\) be signs and
let \(A_1, \dots ,A_q\) be distinct atoms in \((\Omega, \Sigma,\mu)\)
such that \(f_j = \sigma_j \chi_{A_j} / \mu(A_j)\) for all
\(j \leq q\). We can assume without loss of generality that
the \(A_j\)'s are pairwise disjoint.
We define \((x_l^{**})_{l=1}^k\) in \(L_\infty(\mu)\) as follows:
\[ x_l^{**} | A_j = \sigma_j v_{jl}, j \leq q, \text{ while }
   x_l^{**} | \Omega \setminus \cup_{j \leq q} A_j  =
   z_{nl} | \Omega \setminus \cup_{j \leq q} A_j ,\]
where we regard \(z_{nl}\) as an element of \(X^{**}= L_\infty (\mu)\).
Clearly, \(\|x_l^{**}\| \leq 1\) and \(x_l^{**}(f_j)=
\int_{A_j} \sigma_j v_{jl} \sigma_j \chi_{A_j}/ \mu(A_j) d \mu 
=v_{jl}\), for all \(j \leq q\) and \(l \leq k\).

Given signs \(\rho_1 , \dots , \rho_k\) and \(i \leq n\) we have that
\begin{align}
&\biggl \| \sum_{s=1}^k \rho_s \biggl [ \biggl (
\sum_{l=1}^k a_{ils} x_l^{**}
\biggr ) - z_{is} \biggr ] \bigl | \Omega \setminus \cup_{j \leq q} A_j
\biggr \| \leq \label{3} \\
&\biggl \| \sum_{s=1}^k \rho_s \biggl [ \biggl (\sum_{l=1}^k a_{ils} z_{nl}
\biggr ) - z_{is} \biggr ] \biggr \| 
< \delta_i \text{ by \eqref{1} }. \notag
\end{align}
We next fix signs \(\rho_1 , \dots , \rho_k\),
\(i \leq n\) and \(j \leq q\). We set
\[ H_{ij}(t) = \sum_{s=1}^k \rho_s \biggl [ \biggl ( \sum_{l=1}^k a_{ils}
   \sigma_j v_{jl} \biggr ) - z_{is}(t) \biggr ], \,
t \in A_j.\]
We claim that \(|H_{ij}| < \delta_i\), \(\mu\)-almost everywhere
in \(A_j\). Indeed, note first that
\(\int_{A_j} z_{is} d \mu = \sigma_j \mu(A_j) f_j (z_{is})\)
for all \(s \leq k\), and thus
\begin{align} 
\biggl | \int_{A_j} H_{ij} d \mu \biggr | &=
\biggl | \sum_{s=1}^k \rho_s \biggl [ \biggl ( \sum_{l=1}^k a_{ils}
\sigma_j v_{jl} \mu(A_j) \biggr ) - \sigma_j \mu(A_j) f_j (z_{is})
\biggr ] \biggr | \label{4} \\
&\leq \mu(A_j) \sum_{s=1}^k \biggl | \biggl ( \sum_{l=1}^k a_{ils} v_{jl}
\biggr ) - f_j (z_{is}) \biggr | < \delta_i \mu(A_j), 
\text{ by \eqref{2} }. \notag
\end{align}
On the other hand, setting
\(B_{ij}= \{ t \in A_j : \, |H_{ij}(t)| \geq \delta_i\}\)
and taking in account that \(A_j\) is an atom, we infer that
\(\mu(B_{ij})=0\). Indeed, otherwise, \(\mu(B_{ij})=\mu(A_j)\)
and thus \(|H_{ij}| \geq \delta_i\), \(\mu\)-almost everywhere
on \(A_j\). But also, as \(A_j\) is an atom, \(H_{ij}\) has a
constant sign \(\mu\)-almost everywhere on \(A_j\) and so
\(|\int_{A_j} H_{ij} d \mu | \geq \delta_i \mu(A_j)\),
contradicting \eqref{4}. Therefore, \(\mu(B_{ij})=0\)
and hence \(|H_{ij}| < \delta_i\), \(\mu\)-almost everywhere
in \(A_j\), as claimed.
We conclude that 
\(\bigl \| \sum_{s=1}^k \rho_s 
\bigl [ \bigl (\sum_{l=1}^k a_{ils} x_l^{**} \bigr )
 - z_{is} \bigr ] | A_j \bigr \| \leq \delta_i\), for all \( i \leq n\),
\(j \leq q\), and all choices of signs \((\rho_s)_{s=1}^k\).
Combining with \eqref{3} we deduce that
\(\bigl \| \sum_{s=1}^k \rho_s 
\bigl [ \bigl (\sum_{l=1}^k a_{ils} x_l^{**} \bigr )
 - z_{is} \bigr ] \bigr \| \leq \delta_i\), for all \( i \leq n\)
and all choices of signs \((\rho_s)_{s=1}^k\).

We next set \(W= [\{x_l^{**}: l \leq k\} \cup \{ z_{is}:
i \leq n, s \leq k\}]\) and choose \(0 < \delta < \epsilon \).
Theorem \ref{lr} yields a linear operator
\(U \colon W \to X\), \(\|U\| \leq 1 + \delta\),
so that \(U | X \cap W = id_{X \cap W}\) and
\(g(f_j)= f_j (Ug)\), for all \(g \in W\) and \(j \leq q\).
Setting \(x_l = Ux_l^{**} / (1+ \delta)\), \(l \leq k\),
we obtain that for every choice of signs \(\rho_1, \dots , \rho_k\)
and all \(i \leq n\),
\begin{align}
&\biggl \| \sum_{s=1}^k \rho_s 
\biggl [ \biggl (\sum_{l=1}^k a_{ils} x_l \biggr )
 - z_{is} \biggr ] \biggr \| =
\biggl \| \sum_{s=1}^k \rho_s 
\biggl [ \biggl (\sum_{l=1}^k a_{ils} U x_l^{**} /(1 +\delta) \biggr )
 - Uz_{is} \biggr ] \biggr \| \notag \\
&\leq (1 +\delta)
\biggl \| \sum_{s=1}^k \rho_s 
\biggl [ \biggl (\sum_{l=1}^k a_{ils} x_l^{**} /(1 +\delta) \biggr )
 - z_{is} \biggr ] \biggr \| \notag \\
&\leq \biggl \| \sum_{s=1}^k \rho_s 
\biggl [ \biggl (\sum_{l=1}^k a_{ils} x_l^{**}  \biggr )
 - (1 +\delta) z_{is} \biggr ] \biggr \| \notag \\
&\leq \biggl \| \sum_{s=1}^k \rho_s 
\biggl [ \biggl (\sum_{l=1}^k a_{ils} x_l^{**} \biggr )
 - z_{is} \biggr ] \biggr \| + 
\delta \biggl \| \sum_{s=1}^k \rho_s z_{is} \biggr \| 
\leq \delta_i +  \delta, \text{ as }
\|T_i\| \leq 1. \notag
\end{align}
Thus \(\sum_{s=1}^k \bigl | x^* \bigl (
\sum_{l=1}^k a_{ils} x_l \bigr ) - x^*(z_{is}) \bigr |
\leq \delta_i + \delta\), for all \(x^* \in B_{X^*}\).
If we define \(S \colon X^* \to V\) by
\(Sx^* = \sum_{l=1}^k x^*(x_l) e_l\), we see that
\(S\) is \(w^*\)-continuous and \(\|S\| \leq 1\).
Indeed, for the latter assertion we observe that
for every choice of signs \(\rho_1, \dots , \rho_k\),
\(\|\sum_{l=1}^k \rho_l x_l^{**}\| \leq 1\), by
the definition of the sequence \((x_l^{**})_{l=1}^k\)
and the fact that \(\|v_j\| \leq 1\) for \(j \leq q\),
and \(\|T_n\| \leq 1\). It follows now, by the choice
of \(U\), that \(\|\sum_{l=1}^k \rho_l x_l\| \leq 1\)
for every choice of signs \(\rho_1, \dots , \rho_k\),
and therefore \(\|S\| \leq 1\).

We deduce that \(\|R_i Sx^* - T_i x^* \| \leq 
\delta_i + \delta\) for all \(x^* \in B_{X^*}\)
and every \(i \leq n\). Hence \(\|R_i S - T_i \| < 
\delta_i + \epsilon\) for all \(i \leq n\).
Finally, 
\begin{align}
&\|Sf_j - v_j \| = \sum_{l=1}^k | f_j (x_l) - v_{jl}|
\notag \\
&=\sum_{l=1}^k |(1+ \delta)^{-1} f_j (Ux_l^{**}) - v_{jl}| =
\sum_{l=1}^k |(1+ \delta)^{-1} x_l^{**}(f_j) - v_{jl}|
\notag \\
&= \delta (1+ \delta)^{-1} \|v_j\| < \epsilon. \notag
\end{align}
The proof of the proposition is now complete.
\end{proof}
\begin{proof}[{\bf Proof of Theorem \ref{T:2}}.]
We first choose \(\theta_i < \delta_i\) so that
\(\|R_i T_n - T_i \| < \theta_i\) and
\(\|R_i v_j - T_i f_j \| < \theta_i\) for all \(j \leq q\)
and \(i \leq n\). We then choose \(0 < \epsilon_0 <
(\delta_i - \theta_i)/4\) for all \(i \leq n\), and a 
sequence \((\epsilon_m)\) of positive scalars such that
\(\sum_{m=1}^\infty \epsilon_m < \epsilon_0\).

Proposition \ref{P:2} yields a \(w^*\)-continuous linear
operator \(S_0 \colon X^* \to V\), \(\|S_0\| \leq 1\),
such that \(\|R_i S_0 -T_i\| < \theta_i + \epsilon_0\) for
\(i \leq n\), and \(\|S_0 f_j - v_j \| < \epsilon_0\),
for \(j \leq q\).

We next apply Proposition \ref{P:2} for ``\(n\)''\(=1\),
``\(R_1\)''\(= id_V\), ``\(T_1\)''\( =S_0\), 
``\(\delta_1\)''\(=\epsilon_0\)
and ``\(\epsilon\)''=\(\epsilon_1\), to obtain a \(w^*\)-continuous
linear operator \(S_1 \colon X^* \to V\), \(\|S_1\| \leq 1\),
so that \(\|S_0 - S_1\| < \epsilon_0 + \epsilon_1\) and
\(\|S_1 f_j - v_j\| < \epsilon_1 \), for \( j \leq q\).
Continuing in this fashion we construct \(w^*\)-continuous
linear operators \(S_m \colon X^* \to V\) with \(\|S_m\| \leq 1\)
and such that 
\[\|S_{m-1} - S_m \| < \epsilon_{m-1} + \epsilon_m, \,
  \|S_m f_j - v_j \| < \epsilon_m, \, \text{ for all }
  j \leq q, \, m \in \mathbb{N}.\]
Clearly, the sequence of operators \((S_m)\) converges in norm
to a \(w^*\)-continuous linear operator \(T \colon X^* \to V\)
such that \(\|T\| \leq 1\) and \(Tf_j = v_j\), \(j \leq q\).
In addition to that we have
\[\|S_m - S_0\| < \epsilon_0 + 2 \sum_{i=1}^{m-1} \epsilon_i
  + \epsilon_m, \, m > 1,\]
and thus \(\|T - S_0\| < 3 \epsilon_0\). We conclude that
\begin{align}
\|R_i T - T_i \| &\leq \|R_i (T - S_0)\| + \|R_i S_0 - T_i\|
\notag \\
&< \theta_i + 4 \epsilon_0 < \delta_i, \, i \leq n. \notag 
\end{align}
\end{proof}
\begin{Cor} \label{C:4}
Let \(X\) be an \(L_1\)-predual and let \(V\) be
a subspace of \(X^*\) isometric to \(\ell_1^k\).
Assume that there exist
collections of vectors \((f_j)_{j=1}^q \subset 
\mathrm{ext}(B_{X^*})\) and \((v_j)_{j=1}^q \subset B_V\),
with \((f_j)_{j=1}^q \) linearly independent, so that
every linear operator \(R \colon X^* \to X^*\) satisfying
\(Rf_j = v_j\) for all \(j \leq q\), also satisfies
\(R|V = id_V\). Then there exists a 
\(w^*\)-continuous contractive projection \(P\colon X^* \to V\),
such that \(Pf_j = v_j\) for all \(j \leq q\).
\end{Cor}
\begin{proof}
We apply Theorem \ref{T:2} for \(n=1\), \(R_1=id_V\),
\(T_1 = \mathbf{0}\) and \(\delta_1 > 1\), to obtain
a \(w^*\)-continuous linear operator \(P \colon X^* \to V\),
\(\|P\| \leq 1\), so that \(Pf_j = v_j\) for all \(j \leq q\).
Our assumptions yield that \(P\) is the desired projection.
\end{proof}
\begin{remark}
We note that Lemma 3.1 and Corollary 3.2 of \cite{JR} yield
for every \(\epsilon > 0\) a \(w^*\)-continuous projection
\(P \colon X^* \to V\), \(\|P\| \leq 1 + \epsilon \), 
such that \(Pf_j = v_j\) for all \(j \leq q\).
\end{remark}
\section{Main results}
In this section we present the proofs of the results
mentioned in the introduction.
\begin{proof}[{\bf Proof of Theorem \ref{T:1}}.]
Assume \(K\) is infinite and let \((e_n^*)\) be an enumeration of \(K\).
(The argument for finite \(K\) is implicitly contained in the proof
of the infinite case.) It is clear that \((e_n^*)\) is isometrically
equivalent to the usual \(\ell_1\)-basis. Set 
\(Y_n = [e_1^*, \dots e_n^*]\), \(n \in \mathbb{N}\) and
let \((\delta_n)\) be a null sequence of positive scalars.
We shall inductively construct \(w^*\)-continuous contractive
projections \(P_n \colon X^* \to Y_n\) such that
\(\|P_k P_n - P_k\| < \delta_k\), whenever \(k \leq n\).
\(P_1\) is selected by applying Corollary \ref{C:4}
for the subspace \(Y_1\) and the vectors
\(f_1 = v_1 = e_1^*\). Suppose \((P_i)_{i=1}^n\) have been selected
so that \(\|P_i P_j - P_i\| < \delta_i\) whenever
\(i \leq j \leq n\). Apply Theorem \ref{T:2} for 
``\(V\)''\(=Y_{n+1}\), ``\(T_i\)''\(=P_i\),
``\(R_i\)''\(=P_i |Y_{n+1}\), \(i \leq n\), and the collections
of vectors ``\((f_j)_{j=1}^q\)''\(=(e_j^*)_{j=1}^{n+1}\),
``\((v_j)_{j=1}^q\)''\(=(e_j^*)_{j=1}^{n+1}\),
in order to obtain a \(w^*\)-continuous linear operator
\(P_{n+1} \colon X^* \to Y_{n+1}\), \(\|P_{n+1}\| \leq 1\),
such that \(P_{n+1}e_j^*=e_j^*\) for all \(j \leq n+1\),
and \(\|P_i P_{n+1} - P_i \| < \delta_i\) for all \( i \leq n\).
Clearly, \(P_{n+1}\) is the required projection onto \(Y_{n+1}\).
This completes the inductive construction. The assertion
of the theorem now follows from Proposition \ref{P:1}.
\end{proof}
\begin{proof}[{\bf Proof of Corollary \ref{C:1}.}]
Clearly, \(K\) is linearly independent. When \(K\)
is finite the assertion follows immediately from
Theorem \ref{T:1} as \([K]\) is isometric to
\(\ell_1^{|K|}= C(K)^*\). If \(K\) is infinite let
\((e_n^*)\) be an enumeration of \(K\) and set
\(Y=[(e_n^*)]\). Of course \((e_n^*)\) is isometrically
equivalent to the usual \(\ell_1\)-basis, and applying
the Choquet representation and the Krein-Millman theorems,
we infer that \(\overline{\mathrm{co}}^{w^*} (K \cup -K) = B_Y\).
A classical result \cite{KS} yields that \(Y\) is \(w^*\)-closed
in \(X^*\). It is not hard to see (cf. also Lemma 2 of \cite{A2})
that \(Y\) is \(w^*\)-isometric to \(C(K)^*\). The result follows
from Theorem \ref{T:1}.
\end{proof}
\begin{proof}[{\bf Proof of Corollary \ref{C:2}.}]
We regard \(B_{X^*}\) in its \(w^*\)-topology and set
\(H= \mathrm{ext}(B_{X^*})\). Since \(X\) is separable and
\(X^*\) is non-separable, \(H\) is an uncountable 
\(G_{\delta}\)-subset of \(B_{X^*}\). It follows that
\(H\) is an uncountable Polish space in its relative \(w^*\)-topology.
We will show that \(H\) contains a \(w^*\)-compact subset \(L\)
homeomorphic to the Cantor set \(\Delta\), such that
\(L \cap (-L) = \emptyset\). Indeed, let
\(\sigma \colon H \to H\) denote negation (\(\sigma h = -h\)).
Then \(\sigma\) is a fixed-point free homeomorphism on the
uncountable Polish space \(H\) and therefore there exists
an uncountable relatively open subset \(G\) of \(H\) such
that \(G \cap \sigma G = \emptyset\). By a classical result
\(G\) contains a compact subset \(L\) homeomorphic to the Cantor
set which clearly satisfies \(L \cap (-L) = \emptyset\).

Of course \(L\) contains homeomorphs of all countable compact
metric spaces. Corollary \ref{C:1} now yields that \(X\)
contains subspaces isometric to \(C(K)\), for every  
countable compact metric space \(K\). Because \(X\) is separable,
a result of Bourgain (Proposition 9 of \cite{Bo}) implies
that \(X\) contains a subspace isometric to \(C(\Delta)\).
The existence of a contractively complemented subspace of
\(X\) isometric to \(C(\Delta)\) now follows from a result
of Pelczynski \cite{P2}.
\end{proof} 
\begin{proof}[{\bf Proof of Theorem \ref{T:3}}]
We first choose an infinite \(w^*\)-convergent subsequence
\((e_m^*)_{m \in M}\) of \((e_n^*)\) and set
\(x_0^* = w^*-\lim_{m \in M} e_m^*\). Clearly, 
\(\|x_0^*\| \leq 1\). If \(x_0^* = \mathbf{0}\) then
\(Z=[(e_m^*)_{m \in M}]\) is \(w^*\)-closed in \(X^*\)
by Lemma 1 of \cite{A2}. We deduce from Theorem \ref{T:1}
that \(Z\) is the range of a \(w^*\)-continuous contractive
projection in \(X^*\). It is easy to see that if \((r_n)\)
is an enumeration of \(M\) then the subspace
\(Y=[(e_{r_{2n}}^* - e_{r_{2n-1}}^*)]\) is the range of a 
\(w^*\)-continuous contractive projection in \(Z\). 
Hence by composing the projections previously obtained we see
that the assertion of the theorem holds in this case.

We shall next deal with the case of \(x_0^* \ne \mathbf{0}\).
Suppose that \(x_0^* = \sum_{j=1}^\infty a_j e_j^*\)
and choose a sequence of positive scalars \((\epsilon_i)\)
such that \(\sum_{i=1}^\infty \epsilon_i < 1\).
Choose also \(n_1 \in \mathbb{N}\) so that
\(\sum_{j > n_1} |a_j| < \epsilon_1 
\sum_{j \leq n_1} |a_j| \).
We shall inductively construct increasing sequences
\((m_k)_{k=1}^\infty \subset M\) and \((n_k) _{k=1}^\infty \subset
\mathbb{N}\) with \(n_k < m_{2k-1} < m_{2k} < n_{k+1}\),
and \(w^*\)-continuous contractive projections
\(P_k \colon X^* \to Y_k\), where 
\(Y_k = [u_i^*: i \leq k]\) and 
\(u_k^* = (e_{m_{2k}}^* - e_{m_{2k-1}}^*)/2\),
\(k \in \mathbb{N}\), so that the following conditions are
fulfilled:
\begin{align}
&\sum_{j > n_i} |a_j| < \epsilon_i \sum_{j \leq n_1} |a_j|, \,
i \in \mathbb{N}. \label{11} \\
&P_1e_j^* = \mathbf{0}, \, j \leq n_1, \, \text{ while }
\|P_ie_j^*\| < \sum_{l < i} \epsilon_l, \, j \leq n_1, \, i \geq 2.
\label{12} \\
&P_i \biggl ( \sum_{j \leq n_i} a_j e_j^* \biggr )=\mathbf{0}, \,
i \in \mathbb{N}. \label{13} \\
&P_ie_{m_{2j}}^* = u_j^*, \, P_ie_{m_{2j-1}}^* = - u_j^*, \,
j \leq i, \, i \in \mathbb{N}. \label{14} \\
&\|P_i P_j - P_i \| < \sum_{l=i}^{j-1} \epsilon_l, \,
i < j \text{ in } \mathbb{N}. \label{15} \\
&\|P_i e_{m_j}^*\| < \epsilon_i, \, j \in \{2l-1,2l\}, \,
i < l \text{ in } \mathbb{N}. \label{16} 
\end{align}
Once this is accomplished, condition \eqref{15} will enable us
to apply Proposition \ref{P:1} and deduce that
\(Y=[(u_k^*)]\) is the range of a \(w^*\)-continuous projection
in \(X^*\). Note that \(Y\) is \(w^*\)-closed in \(X^*\) by
Lemma 1 of \cite{A2} as \((u_k^*)\) is isometrically equivalent 
to the usual \(\ell_1\)-basis.

We first choose \(m_1 < m_2\) in \(M\) with \(m_1 > n_1\),
and apply Corollary \ref{C:4} for ``\(V\)''\(=Y_1\),
``\(q\)''\(=n_1 +2\), ``\(f_j\)''\( = e_j^*\) (\(j \leq n_1\)),
``\(f_{q-1}\)''\(=e_{m_1}^*\), ``\(f_q\)''\(=e_{m_2}^*\), and
``\(v_j\)''\(=\mathbf{0}\) (\(j \leq n_1\)),
``\(v_{q-1}\)''\(=-u_1^*\), ``\(v_q\)''\(=u_1^*\).
We obtain a \(w^*\)-continuous contractive projection
\(P_1 \colon X^* \to Y_1\) such that \(P_1e_j^* =\mathbf{0}\)
for all \(j \leq n_1\).

Suppose that we have constructed \((n_i)_{i=1}^k\),
\((m_i)_{i=1}^{2k}\) and \((P_i)_{i=1}^k\) so that
conditions \eqref{11}-\eqref{16} are satisfied.
We next choose \(n_{k+1} > m_{2k}\) so that \eqref{11}
is satisfied for \(i=k+1\). By \eqref{11} and \eqref{13}
of the induction hypothesis we infer that
\(\|P_i x_0^*\| < \epsilon_i\), for \(i \leq k\).
We can therefore choose \(m_{2k+1} < m_{2k+2}\)
in \(M\) with \(m_{2k+1} > n_{k+1}\), such that
\(\|P_i e_{m_{2k+1}}^*\| < \epsilon_i\) and
\(\|P_i e_{m_{2k+2}}^*\| < \epsilon_i\) for all \(i \leq k\).
Hence \eqref{16} is satisfied for \(l=k+1\).

We next put \(q=n_{k+1} + 2\) and set \(f_j = e_j^*\), for
\(j \leq n_{k+1}\), \(f_{q-1}= e_{m_{2k+1}}^*\) and
\(f_q=e_{m_{2k+2}}^*\). 
We claim that there exist vectors \((v_j)_{j=1}^q\) in
\(B_{Y_{k+1}}\) so that
\begin{align}
&\|v_j\| < \sum_{l=1}^k \epsilon_l, \, j \leq n_1.
\label{17} \\
&v_j = P_k e_j^*, \, j \in (n_1 , n_{k+1}], \,
 v_{q-1}= - u_{k+1}^*, \, v_q = u_{k+1}^*.
\label{18} \\
&\sum_{j \leq n_{k+1}} a_j v_j = \mathbf{0}. \label{19} \\
&\|P_i f_j - P_i v_j \| < \sum_{l=i}^k \epsilon_l, \,
i \leq k, \, j \leq q. \label{20}
\end{align} 
Having achieved this and taking in account
\eqref{15} of the induction hypothesis,
we employ Theorem \ref{T:2}
for ``\(V\)''\(=Y_{k+1}\), ``\(n\)''\(=k\),
``\(\delta_i\)''\(=\sum_{l=i}^k \epsilon_l\),
``\(T_i\)''\( = P_i\), ``\(R_i\)''\(=P_i |Y_{k+1}\)
(\(i \leq k\)), and the collections of vectors
\((f_j)_{j=1}^q\), \((v_j)_{j=1}^q\) described above,
to find a \(w^*\)-continuous linear operator
\(P_{k+1} \colon X^* \to Y_{k+1}\), \(\|P_{k+1}\| \leq 1\),
such that \(\|P_i P_{k+1} - P_i\| < \sum_{l=i}^k \epsilon_l\),
for all \(i \leq k\), and \(P_{k+1}f_j = v_j\) for all
\(j \leq q\). It is easy to verify that \(P_{k+1}\)
is a projection onto \(Y_{k+1}\) so that
\((n_i)_{i=1}^{k+1}\),
\((m_i)_{i=1}^{2k+2}\) and \((P_i)_{i=1}^{k+1}\) satisfy
conditions \eqref{11}-\eqref{16}. 

The collection \((v_j)_{j=n_1 + 1}^q\) is explicitly defined 
in \eqref{18}. It remains to define \((v_j)_{j=1}^{n_1}\).
We first choose scalars \((b_{il})\), where \(i \leq k\)
and \(l \in (n_k , n_{k+1}]\), such that
\(P_k e_l^* = \sum_{i=1}^k b_{il} u_i^*\).
Note that \(\sum_{i=1}^k |b_{il}| \leq 1\) for every
\(l \in (n_k , n_{k+1}]\) since \(\|P_k\|=1\).
We also define scalars \((\rho_i)_{i=1}^k\) by
\[ \rho_i = \biggl (- \sum_{l \in (n_k , n_{k+1}]}
   a_l b_{il} \biggr ) \biggl / \sum_{j \leq n_1} |a_j|,\]
and set \(\rho_{ij} = sg(a_j) \rho_i\) for \(i \leq k\)
and \(j \leq n_1\). Observe that
\(\sum_{i=1}^k |\rho_i| < \epsilon_k\), by \eqref{11}
of the induction hypothesis. We now set
\[ v_j = P_k e_j^* + \sum_{i=1}^k \rho_{ij} u_i^*, \,
j \leq n_1.\]
It follows now by \eqref{12} of the induction hypothesis
that \eqref{17} is satisfied.
To establish \eqref{19} we have
\begin{align}
\sum_{j \leq n_{k+1}} a_j v_j &= \sum_{j \leq n_1} a_j v_j +
\sum_{j \in (n_1 , n_{k+1}]} a_j P_k e_j^* 
\notag \\
&=\sum_{j \leq n_1} a_j \sum_{i=1}^k \rho_{ij} u_i^* +
\sum_{j \leq n_{k+1}} a_j P_k e_j^* 
\notag \\
&=\sum_{j \leq n_1}|a_j| \sum_{i=1}^k \rho_i u_i^* +
\sum_{j \leq n_k} a_j P_k e_j^* +
\sum_{j \in (n_k , n_{k+1}]} a_j P_k e_j^* 
\notag \\
&= - \sum_{i=1}^k \sum_{l \in (n_k , n_{k+1}]}
a_l b_{il} u_i^* +
\sum_{j \in (n_k , n_{k+1}]} a_j P_k e_j^* ,
\text{ by } \eqref{13}
\notag \\
&= - \sum_{l \in (n_k , n_{k+1}]}
a_l \sum_{i=1}^k b_{il} u_i^* +
\sum_{j \in (n_k , n_{k+1}]} a_j P_k e_j^* 
\notag \\
&= - \sum_{l \in (n_k , n_{k+1}]} a_l P_ke_l^* +
\sum_{j \in (n_k , n_{k+1}]} a_j P_ke_j^* = \mathbf{0}. \notag
\end{align}
Finally, we show that \eqref{20} holds.
Indeed, when \(j \in \{q-1,q\}\), this is a 
consequence of the choice of \(m_{2k+1}\) and \(m_{2k+2}\).
When \(j \in (n_1, n_{k+1}]\) the assertion follows
from \eqref{15} of the induction hypothesis. 
When \(j \leq n_1\) it follows from \eqref{15}
of the induction hypothesis and the fact that
\(\sum_{i=1}^k |\rho_i| < \epsilon_k\).
\end{proof}
\begin{proof}[{\bf Proof of Corollary \ref{C:3}}]
If \(X^*\) is non-separable the assertion follows from
Corollary \ref{C:2}. If \(X^*\) is separable, then
\(\mathrm{ext}(B_{X^*})\) is countable and \(X^*\)
is isometric to \(\ell_1\). Let \((e_n^*)\)
be a basis for \(X^*\) isometrically equivalent to
the usual \(\ell_1\)-basis, and choose a \(w^*\)-convergent
subsequence \((e_{m_n}^*)\) of \((e_n^*)\) according
to Theorem \ref{T:3}. Let \(Y= [(e_{m_{2n}}^* -
e_{m_{2n-1}}^*)]\). Then it is easy to see 
that \(Y\) is \(w^*\)-isometric to \(c_0^*\).
The result follows from Theorem \ref{T:3}.
\end{proof} 
Our last corollary is a special case of the structural
result for separable \(L_1\)-preduals established in
\cite{MP} and \cite{Ll}.
\begin{Cor}
Suppose that \(X^*\) is isometric to \(\ell_1\).
Then there exists a sequence \((E_n)\) of finite-dimensional
subspaces of \(X\) such that \(E_n \subset E_{n+1}\)
for all \(n\), each \(E_n\) is isometric to 
\(\ell_{\infty}^n\), and \(\cup_{n=1}^\infty E_n\) is dense
in \(X\).
\end{Cor}
\begin{proof}
Let \((e_n^*)\) be a basis for \(X^*\) isometrically
equivalent to the usual \(\ell_1\)-basis. Let
\(Y_n = [e_i^*: i \leq n]\), \(n \in \mathbb{N}\),
and let \((\delta_n)\) be a sequence of positive scalars
such that \(\sum_{n=1}^\infty \delta_n < \infty\). The argument
in the proof of Theorem \ref{T:1} now yields a sequence
\((P_n)\) of \(w^*\)-continuous contractive projections
in \(X^*\) with \(\mathrm{Im} P_n = Y_n\) and such that
\(\|P_k P_n - P_k \| < \delta_k\) whenever \(k \leq n\).
Given \( k \leq n \) in \(\mathbb{N}\), we set
\(Q_k^n = P_k \cdots P_n\). Clearly \(Q_k^n\)
is a \(w^*\)-continuous contractive projection onto \(Y_k\).
Moreover, our assumptions yield that 
\(\|Q_k^n - Q_k^{n+1}\| < \delta_n\), for all \(n \geq k\)
and thus the sequence of operators \((Q_k^n)_{n \geq k}\)
converges in norm to a \(w^*\)-continuous contractive
projection \(Q_k\) from \(X^*\) onto \(Y_k\).
It is easily seen that \(Q_k^m Q_l^n = Q_k^n\), whenever
\(k \leq l \leq m \leq n\) and hence \(Q_k Q_l = Q_k\)
when \(k \leq l\). 

We now let \(E_n = Q_n^* Y_n^*\). \(E_n\) is naturally
identified to a subspace of \(X\) as \(Q_n\) is \(w^*\)-continuous,
and of course it is isometric to \(\ell_\infty^n\) for all
\(n \in \mathbb{N}\). Since \(Q_n Q_{n+1} = Q_n\) we deduce
that \(E_n \subset E_{n+1}\) for all \(n \in \mathbb{N}\).
It is also easily verified that \(Q_n^*\) acts as a contractive
projection from \(X\) onto \(E_n\). Rainwater's theorem now
yields that \(\lim_n Q_n^* x = x\), weakly, for all \(x \in X\)
and thus
\(\cup_{n=1}^\infty E_n\) is dense in \(X\).
\end{proof}
\begin{remark}
\begin{enumerate}
\item We note that the proof of Corollary \ref{C:3}
that appears in \cite{Z} makes use of the structural result
established in \cite{MP} and \cite{Ll}. 
\item Theorem \ref{T:3} is no longer valid if we
consider isomorphic \(\ell_1\)-preduals. Indeed,
it is shown in \cite{BD} that there exist 
isomorphic \(\ell_1\)-preduals which do not contain isomorphic
copies of \(c_0\).
\item According to a result of Fonf \cite{f}, every
Banach space \(X\) such that \(\mathrm{ext}(B_{X^*})\)
is countable contains a subspace isomorphic to \(c_0\). 
\end{enumerate}
\end{remark}
It was shown in \cite{JZ} that every separable \(L_1\)-predual
is isometric to a quotient of \(C(\Delta)\). It is an open
problem \cite{A2} whether every \(\ell_1\)-predual is
isomorphic, or even almost isometric, to a quotient of
\(C(K)\) for some countable compact metric space \(K\).

{\em Question:} Suppose \(X\) is an \(\ell_1\)-predual
such that for some \(\epsilon > 0\) 
and some countable ordinal \(\alpha\),
the \(\epsilon\)-Szlenk
index of \(X\) \cite{S} exceeds \(\omega^\alpha\).
Does \(X\) contain a contractively complemented subspace
isometric to \(C_0(\omega^{\omega^\alpha})\)?
Does \(X\) contain a subspace isomorphic to
\(C(\omega^{\omega^\alpha})\)?

\end{document}